\documentclass[reqno]{amsart}
\usepackage{amssymb}  
\usepackage[T2A]{fontenc}
\usepackage[utf8]{inputenc}
\usepackage[english]{babel} 
\usepackage{cite}     
\usepackage[dvipsnames]{xcolor} % for link color names
\usepackage{hyperref} % clickable URL
\usepackage{breakurl}
\hypersetup{
    colorlinks=true,
    linkcolor=OliveGreen,  
    urlcolor=BlueViolet,
    citecolor=MidnightBlue
    }
\usepackage{multirow} % for multirows in tables
\usepackage{threeparttable} % for footnotes in tables
\usepackage{caption} % to widen table caption

\newtheorem{thm}{Theorem}

\begin{document}

\title[On generation by triality automorphisms]{On generation by triality automorphisms}
\author{Danila O. Revin}%
\address{Danila O. Revin
\newline\indent Sobolev Institute of Mathematics,
\newline\indent 4, Koptyug av.
\newline\indent 630090, Novosibirsk, Russia
} 
\email{revin@math.nsc.ru
\newline\indent\normalfont ORCID: \href{https:\\orcid.org/0000-0002-8601-0706}{0000-0002-8601-0706}}

\author{Andrei V. Zavarnitsine}%
\address{Andrei V. Zavarnitsine
\newline\indent Sobolev Institute of Mathematics,
\newline\indent 4, Koptyug av.
\newline\indent 630090, Novosibirsk, Russia
} 
\email{zav@math.nsc.ru
\newline\indent\normalfont ORCID: \href{https:\\orcid.org/0000-0003-1983-3304}{0000-0003-1983-3304}
}
\thanks{This research was carried out within the State Contract of the Sobolev Institute of Mathematics (FWNF-2026-0017).}
\maketitle
\begin{quote}
\noindent{\sc Abstract. } 
We clarify the structure of subgroups generated by conjugate graph automorphisms of order $3$ of $O_8^+(2)$ and $O_8^+(3)$. As a result, we obtain a correction to a paper by S. Guest which, in turn, plays an important role in proving the solvable analogue of the Baer--Suzuki theorem.

\smallskip



\end{quote}

\section{Introduction}

In 2009--2010, two independent collectives of mathematicians, namely, P.\,Flavell, R.\,Gu\-ral\-nick, and S.\,Guest \cite{FGG}, as well as L.\,Gordeev, F.\,Grunewald, B.\,Kunyavski\u{\i}, and E.\,Plotkin \cite{GGKP,GGKP2} used the classification of finite simple groups to prove a solvable analogue of the Baer--Suzuki theorem. This theorem asserts that an element $x$ of a finite group~$G$ lies in the solvable radical of $G$ if the subgroup $\langle x,x^{g_1},x^{g_2},x^{g_3}\rangle$ is solvable for every ${g_1,g_2,g_3\in G}$. Either collective reduced the proof to the situation where $G$ is almost simple, i.\,e.,\footnote{As usual, we identify the isomorphic groups $S$ and $\operatorname{Inn}S\leqslant \operatorname{Aut}S$.} $S\leqslant G\leqslant\operatorname{Aut}S$ for a finite simple nonabelian group $S$, an automorphism $x$ of $S$ of prime order, and $G=\langle S,x\rangle$. 

An important component of \cite{FGG} is Theorem A$\mbox{}^*$ by S.\,Guest \cite{Guest} which asserts that if~$x$ is an automorphism of odd prime order of $S$ as above then, up to a reasonable number of exceptions, there is $g\in G$ such that the group $\langle x,x^g\rangle$ is not solvable. In the same Theorem A$\mbox{}^*$, a classification of all possible exceptions was proposed. It turns out that, in all exceptional cases, the subgroup $\langle x,x^{g_1},x^{g_2},x^{g_3}\rangle$ is nonsolvable for some  $g_1,g_2,g_3\in S$. Considering what was said in the previous paragraph, Theorem A$\mbox{}^*$ implies the validity of the solvable radical theorem for elements $x$ of odd order. Furthermore, Theorem A$\mbox{}^*$ was used in \cite[Theorem 2.4 and Lemma 2.5]{FGG} to prove the solvable radical theorem also in the case where $x$ is an involutive automorphism\footnote{Subsequently \cite[Theorem~1.6]{Guest2}, Guest also classified all the cases where $x$~is an involutive automorphism of a simple group $S$ and the subgroup $\langle x, x^{g_1}, x^{g_2}\rangle\leqslant\langle S,x\rangle$ is solvable for every $g_1, g_2 \in S$. (Of course, if $x$ is an involution, it generates a solvable subgroup together with any of its conjugates.)}. Guest's Theorem~A$\mbox{}^*$ has many other applications, see, e.\,g., \cite{R11,WGR,RZ24,TV24,AGS,GMT,YRV}.

In this paper, we point out two exceptional cases missing from \cite[Theorem~A$\mbox{}^*$]{Guest}. These are the cases where\footnote{Here and henceforth, we use the standard notation from the ATLAS of Finite Groups \cite{atlas}. In particular, we denote by $O_{2n}^\varepsilon(q)$ the group $\mathrm{P}\Omega_{2n}^\varepsilon(q)$, where $\varepsilon\in\{+,-\}$.} $S=O_8^+(q)$ for $q=2,3$ and $x=\rho$, a triality automorphism, see Section~\ref{s:triality}: 

\begin{thm}\label{cor}
Let $G=\langle S,\rho \rangle$, where  $S= O_8^+(2)$ or $O_8^+(3)$ and $\rho$ is a  triality automorphism of $S$. Then $\langle\rho,\rho^g\rangle$ is a solvable group for every $g\in G$.
\end{thm}

The reason why these cases were omitted from \cite{Guest} is as follows. In \cite[proof of Lemma 7, p. 5917]{Guest}, it is observed that a graph or field-graph automorphism $x$ of order $3$ of $D_4(q)\cong O_8^+(q)$ always normalises but does not centralise a subgroup isomorphic to $G_2(q)$ and thereby induces an automorphism of order~$3$ on $G_2(q)'$. Based on this, it is concluded that $(x,D_4(q))$ cannot be a minimal counterexample to Theorem A$\mbox{}^*$. However, it is not taken into account that, in the case of a graph\footnote{In the case of a field-graph automorphism, $q$ must be a cube, so no new exceptions arise.} automorphism $x$, it is possible for  $G_2(q)'$ with $q=2,3$ to admit an automorphism of order $3$ that appears on the list of exceptions \cite[Table~1]{Guest} (cf. Table~\ref{A**}), and the possibility that such an automorphism of $G_2(q)'$ might be induced by~$x$ is overlooked. 

Because the proof of \cite[Theorem A$\mbox{}^*$]{Guest} is performed by shattering a minimal counterexample, new exceptions to the statement of the theorem, unaccounted for in \cite{Guest}, might {\em a~priori} lead to the emergence of entire new series of such exceptions. Fortunately, the case of a graph automorphism of order $3$ of $O_8^+(q)$, being too specific, does not arise in any subsequent argument of Guest \cite{Guest}.
Therefore, to bridge the gap in \cite[Theorem A$\mbox{}^*$]{Guest}, it is sufficient to find all exceptions that emerge from such automorphisms for $q=2,3$ and  correct accordingly the statement of the theorem. 

A graph automorphism of order~$3$ of $O_8^+(q)$ is not always a triality, but always shares with a triality the same coset of $O_8^+(q)$ in $\operatorname{Aut} O_8^+(q)$, see Section~\ref{s:triality} for details. We will show that every triality automorphism gives rise to exceptional cases in \cite[Theorem~A$\mbox{}^*$]{Guest}, while the nontriality graph automorphisms of order~$3$ do not. This is accomplished by finding the value $\alpha(x)$ for all graph automorphisms~$x$ of order $3$. Here, $\alpha$ is a parameter introduced by R.\,Guralnick and J.\,Saxl \cite{GS} which means the following. Given a nonabelian simple group $S$ and its nonidentity automorphism $x$, the value $\alpha(x)=\alpha^{\phantom{S}}_S(x)$ is the minimum number of conjugates of~$x$ in $\langle S, x\rangle$ that generate a subgroup containing~$S$. This parameter has been broadly studied, e.\,g., in \cite{BM,BM2,DMPZ,FMOBW, RZ24,WGR}, and found numerous applications including those in \cite{FGG,Guest,GGKP,GGKP2}. Clearly, if  $\alpha(x)=2$ then, for some $g\in S$, the subgroup $\langle x, x^g\rangle$ is nonsolvable. 

Recall also that, given a set of primes~$\pi$, a {\em $\pi$-group} is a finite group whose order has all its prime divisors belonging to~$\pi$. Our principal result is as follows.

\begin{thm}\label{main}
Let $S$ be either $O_8^+(2)$ or $O_8^+(3)$ and let $\rho$ be a triality automorphism of $S$ of order $3$. Then 

\begin{enumerate}
\item[$(i)$] For every $g\in S$, the subgroup $\langle \rho, \rho^g\rangle$ of $\langle S,\rho \rangle$ is a $\{2,3\}$-group. Moreover, $\alpha(\rho)=3$. 

\item[$(ii)$] If $\varrho \in S\rho$ of order $3$ is not conjugate to $\rho$ then $\alpha(\varrho)=2$. 
\end{enumerate}
\end{thm}

Observe that, for a graph automorphism $x\in \operatorname{Aut}S$ of order~$3$ of $S= O_8^+(q)$, only the estimate $\alpha(x)\leqslant 6$ is known \cite[p.~541]{GS} in the case of an arbitrary~$q$.

Our proof of Theorem~\ref{main} uses some computations in the computer algebra system \texttt{GAP}~\cite{GAP}. The programming code with comments that accompanies these computations is openly available for verification from \cite{RZGenTr}.

Theorem~\ref{cor} follows\footnote{In Table~\ref{rrgen}, we list the structure of all subgroups of the form $\langle\rho,\rho^g\rangle$ for $g\in S$.} from item~$(i)$ of Theorem~\ref{main} in view of Burnside's theorem\footnote{See \cite[Theorem~III]{Burnside} and \cite[Theorem (3.10)]{Isaacs}.} which asserts that a $\pi$-group is solvable whenever $\pi=\{p,q\}$.

The above remarks, therefore, may be considered sufficient to view \cite[Theorem~A$\mbox{}^*$]{Guest} as proven in its modified form as follows.

\begin{thm}[Guest's Theorem A$\mbox{}^*$ corrected]\label{GuestA*}
Let $G$ be a finite almost simple group with socle $S$. Suppose that
$x$ is an element of odd prime order in $G$. Then one of the following holds.
\begin{itemize}
    \item[$(i)$] There exists $g\in G$ such that $\langle x, x^g\rangle$ is not solvable.
     \item[$(ii)$] $x^3 = 1$ and $(x, S)$ belongs to a short list of exceptions given in Table~{\rm \ref{A**}}.
Moreover, there exist $g_1, g_2 \in G$ such that $\langle x, x^{g_1} , x^{g_2}\rangle$ is not solvable, unless
$S\cong  U_n(2)$ or $S\cong  S_{2n}(3)$. In any case, there exist $g_1, g_2, g_3 \in G$ such that $\langle x, x^{g_1} , x^{g_2} , x^{g_3}\rangle$ is not solvable.
\end{itemize}
\end{thm}
\begin{table}[htb]
\caption{List of exceptions to Theorem A*.\label{A**}}
\begin{tabular}{|c|c|}
\hline
$S\vphantom{S^{S^S}}$ & $x$\\
\hline
$L_n(3)$, $n > 2\vphantom{S^{S^S}}$ & transvection\\
$S_{2n}(3)$, $n > 1$ & transvection\\
$U_n(3)$,  $n > 2$ & transvection\\
$U_n(2)$, $n > 3$ & pseudoreflection of order $3$\\
$O^\varepsilon_n(3)$, $n > 6$,  $\varepsilon\in\{+,-,\circ\}$ &  long root element\\

$E_6(3)$, $E_7(3)$, $E_8(3)$, $F_4(3)$, ${}^3D_4(3)$, ${}^2E_6(3)$ &  long root element \\
$G_2(3)$ &  long or short root element\\
$G_2(2)'\cong U_3(3)$ & transvection\\
$O^+_8(2)$, $O^+_8(3)\vphantom{g_{g_{g_g}}}$ &  triality\\
\hline
\end{tabular}
\end{table}

Finally, we remark that the discovery of new exceptions to \cite[Theorem A$\mbox{}^*$]{Guest} by no means contradicts the main result of the paper \cite{FGG} and leaves the proof of its key assertions valid. The results in \cite{R11,WGR,RZ24,YRV}, where an analogue of the  Baer--Suzuki theorem for the $\pi$-radical is studied and where \cite[Theorem A$\mbox{}^*$]{Guest} is significantly used, remain valid as well for the same reason. We also note that the values of $\alpha$ obtained in Theorem~\ref{main} will find independent application in the study of this analogue of the Baer–Suzuki theorem \cite[Conjectures 1,\,2]{YRV}.

\section{Triality and other graph automorphisms of order~3}\label{s:triality}

We fix a finite field~$\mathbb{F}_q$ with $q$ elements. Throughout this section, 
$S$ will stand for a simple orthogonal group $O_8^+(q)$ or, in the Lie notation, its isomorphic Chevalley group $D_4(q)$ over~$\mathbb{F}_q$, see \cite[Section~4.4]{Carter}. In this section, we will recall some facts about graph automorphisms of order $3$ of $S$ and make clearer the terminology we use.

$S$ has a rather unusual for simple groups group of outer isomorphisms due to the existence of a so-called triality automorphism. There is some inconsistency in the literature regarding what is meant by the term ``triality''.\footnote{
See, e.\,g., footnote \ref{klfoot} on the use of the term ``triality'' in \cite{Klei}. We also mention the so-called {\em groups with triality} which have been widely studied due to their close relation to Moufang loops \cite{Glaub,GZ,Liebeck}. $O_8^+(q)$ provides an important example of such groups. The presence of triality for $O_8^+(q)$ correlates with the existence of a series of finite nonassociative simple Moufang loops. The notation $\rho$ which we use for triality was borrowed from the papers on groups with triality.}
We will try to keep our terminology close to that of~\cite{Guest} and, in particular, use ``triality'' in the following sense. The Dynkin diagram of a Lie algebra of type $D_4$ has a symmetry of order $3$ which, according to \cite[Proposition~12.2.3]{Carter}, can be canonically lifted to an automorphism~$\rho$ of order $3$ of the simple group $D_4(q)\cong S$. This automorphism~$\rho$, along with any of its conjugates in $\operatorname{Aut}S$, will be called \emph{triality}.\footnote{See also \cite[Theorem 2.5.1(d1)]{GLS}. An alternative approach to construct triality can be found in \cite[Section~4.7.2]{Wilson}.} Guest \cite[p.~5917]{Guest} calls such automorphisms ``standard triality''.

With respect to the conjugation action of $S$, the set of all trialities decomposes into two orbits 
for $q$ even and eight orbits for $q$ odd. Namely, there is a one-to-one correspondence between the $S$-conjugacy classes of trialities and the $3$-cycles in the image of the epimorphism 
$$ \mu:\operatorname{Aut}S\to\left\{
\begin{array}{ll}
   S_3  & \text{ if } q \text{ is even,}\\
    S_4 & \text{ if } q \text{ is odd. }
\end{array}
\right.$$
Here, $S_3$ (also as a subgroup of $S_4$) is the full group of symmetries of the Dynkin diagram of type $D_4$, and the normal Klein $4$-subgroup of $S_4$ (which is complemented by  $S_3$) is the image under $\mu$ of the group $\operatorname{Inndiag}S$ of inner-diagonal automorphisms.

The trialities are mapped by $\mu$ to $3$-cycles.\footnote{As a caveat, we note that, say, P.\,Kleidman \cite{Klei} refers to any automorphism with this  property as ``triality'' regardless of whether it has order~$3$.\label{klfoot}}
However, not only the trialities have this property, but also, for example, the nontriality {\em graph automorphisms of order~$3$} which play an important role in this paper. An automorphism of order~$3$ of~$S$ is\footnote{This definition agrees with \cite[Definition~2.5.13(b3)]{GLS}, as does the paper \cite{Guest}.} a graph automorphism if it lies in the same coset of $S$ in
$\operatorname{Aut}S$ with some triality. Besides  $\rho^S$, there is 
precisely one\footnote{See \cite[Table~4.7.3A  for ${3\nmid q}$ and Proposition~4.9.2  for ${3\mid q}$\,]{GLS}.} $S$-conjugacy class of nontriality elements of order $3$ in the coset $S\rho$ whose representative will be denoted by~$\varrho$. Moreover, we have\footnote{{\em Ibid.} and \cite[Proposition~1.4.1]{Klei}.}   
\begin{equation}\label{crv}
C_S(\rho)\cong G_2(q)\quad \text{and} \quad 
C_S(\varrho)\cong\left\{
\begin{array}{rl}
\mathrm{PGL}_3(q) & \text{ if } q\phantom{-}\equiv \phantom{-}1 \pmod 3,\\
\mathrm{PGU}_3(q) & \text{ if } q\phantom{-}\equiv {-}1 \pmod 3,\\ 
\left[q^5\right]\!\!.\mathrm{SL}_2(q) & \text{ if } q\phantom{-}\equiv \phantom{-}0 \pmod 3.
\end{array}
\right.
\end{equation}
 In any case, $|C_S(\rho)|\ne |C_S(\varrho)|$.    

\section{Proof of main result}
In this section, we prove Theorem \ref{main} keeping the above notation.
\begin{proof} 
First, let $S\cong O_8^+(2)$. The extension $\langle S,\rho\rangle \cong {O_8^+(2)\!:\!3}$ 
has no available representation in the online atlas \cite{AtlRep}, while \texttt{GAP} \cite{GAP} includes it as a primitive permutation group of degree $1575$. Theoretically, $\langle S,\rho\rangle$ has a matrix $24$-dimensional representation over $\mathbb{F}_2$ which may be of independent interest, and so we briefly outline its explicit construction.\footnote{More generally, we provide in \cite{RZGenTr} generating $24\!\times\! 24$-matrices over $\mathbb{F}_2$ for the extension of $O_8^+(2)$ by the full group $S_3$ of graph automorphisms.} 

\begin{table}[htb]
\caption{Generators $x$, $y$ of $O_8^+(2)$ and their images   
under the triality $\rho$.\label{gensO8}}
\begin{tabular}{r@{}lr@{}l}
$x=$&
$\left({\scriptstyle
\begin{array}{cccccccc}
    1&.&1&1&1&1&.&1\\
     .&.&1&.&1&1&.&1\\
     .&.&1&.&1&.&.&1\\
     .&.&.&.&1&.&.&1\\
     .&.&.&.&1&.&.&.\\
     .&1&.&.&1&.&.&.\\
     .&1&.&.&.&.&1&.\\
     .&.&.&1&.&.&1&.
\end{array}
}\right)$&
$y=$&
$\left({\scriptstyle
\begin{array}{cccccccc}
    1&1&.&.&.&.&.&.\\
     .&1&.&.&.&.&.&.\\
     .&.&1&.&.&.&.&.\\
     .&.&.&1&.&.&.&.\\
     .&.&.&1&1&.&.&.\\
     .&.&.&.&1&1&.&.\\
     .&.&.&.&.&.&1&.\\
     1&1&.&.&.&.&.&1
\end{array}
}\right)$\\[50pt]
$x^\rho=$&
$\left({\scriptstyle
\begin{array}{cccccccc}
    .&.&1&.&.&.&.&.\\
     1&1&.&1&.&.&.&.\\
     1&1&.&.&.&.&.&.\\
     .&1&.&.&.&.&.&.\\
     1&1&.&1&.&.&1&.\\
     1&.&1&.&.&.&.&1\\
     1&.&.&1&1&.&.&1\\
     .&.&.&1&1&1&.&.
\end{array}
}\right)$&
$y^\rho=$&
$\left({\scriptstyle
\begin{array}{cccccccc}
    1&.&.&.&.&.&.&.\\
     .&1&.&.&.&.&1&.\\
     1&.&1&.&.&1&.&.\\
     .&.&.&1&.&.&.&.\\
     .&1&.&.&1&.&.&.\\
     1&.&.&.&.&1&.&.\\
     .&.&.&.&.&.&1&.\\
     .&.&.&.&.&.&.&1
\end{array}
}\right)$
\end{tabular}
\end{table}

Since $S$ is isomorphic to $\Omega_8^+(2)$, it has a natural $8$-dimensional representation over $\mathbb{F}_2$. We check in~\cite{RZGenTr} that the $8\times 8$-matrices $x$ and~$y$ from Table~\ref{gensO8}, where the dots stand for zeros, generate a group isomorphic to~$S$. Now, $\langle S,\rho\rangle$, being an overgroup of $S$ of index $3$, has a $24$-dimensional induced representation. In order to construct it explicitly, we define $\rho$ by its action on the generators 
\begin{equation}\label{xyrho}
x\to x^\rho, \qquad y\to y^\rho,
\end{equation}
where $x^\rho$ and $y^\rho$ are also given in Table \ref{gensO8}. More precisely, we check that $x^\rho,y^\rho\in S= \langle x,y \rangle$,  the map (\ref{xyrho}) indeed extends to an automorphism of $S$ using \texttt{GAP}'s 
function \texttt{GroupHomomorphismByImages()}, $\rho$ has order~$3$, and that 
$$|C_S(\rho)|=12096=|G_2(2)|.$$ 
According to \cite{atlas}, $S$ has no elements of order $3$ with centraliser of such order; hence, (\ref{crv}) implies that $\rho$ is indeed a triality.

The group generated by the block-diagonal $24\times 24$-matrices
\begin{equation}\label{xyind}
\left(
\begin{array}{ccc}
     x&.&.\\
     .&x^{\rho^2}&.\\
     .&.&x^\rho
\end{array}
\right),\qquad 
\left(
\begin{array}{ccc}
     y&.&.\\
     .&y^{\rho^2}&.\\
     .&.&y^\rho
\end{array}
\right)
\end{equation}
is clearly isomorphic to $S$, and its automorphism $\rho$ can be naturally identified with the permutation matrix 
\begin{equation}\label{rhoind}
\left(
\begin{array}{ccc}
     .&I&.\\
     .&.&I\\
     I&.&.
\end{array}
\right),
\end{equation}
where $I$ is the identity $8\times 8$-matrix. Using \texttt{GAP}'s function 

\smallskip\noindent
\texttt{StructureDescription()} 

\smallskip\noindent
we check that the group generated by (\ref{xyind}) and (\ref{rhoind}) has the required structure ${O_8^+(2)\!:\!3}$. This construction allows us to verify that $\langle \rho, \rho^x, \rho^{y}\rangle = \langle S,\rho\rangle$ which yields ${\alpha(\rho)\leqslant 3}$.

To describe the structure of subgroups $\langle \rho,\rho^{g} \rangle$,  $g\in S$, it suffices to assume that $\rho^{g}$ is an orbit representative for $C_S(\rho)$ acting on the conjugacy class $\rho^S$. There are $17$ such orbits, and the corresponding isomorphism types of $\langle \rho,\rho^{g} \rangle$ are listed in Column I of Table \ref{rrgen}. Column II gives the \texttt{GAP} ID for the subgroups (which allows one to construct them uniquely using \texttt{GAP}'s function \texttt{SmallGroup()}) except for those of size $3456$. In Column IV, which addresses both cases $S\cong O_8^+(2)$ and $O_8^+(3)$, we list the sizes of $C_S(\rho)$-orbits on $\rho^S$ whose representatives, together with $\rho$,  generate a subgroup of the given isomorphism type. The symbol ``($\times 3$)'' next to a number means that there are three orbits of this 
size, see \cite{RZGenTr} for details.

\begin{table}[htb]
\begin{threeparttable}
\caption{Subgroups $\langle \rho,\rho^g \rangle$ of  $\langle S,\rho\rangle$ for $S\cong O_8^+(2)$ or $O_8^+(3)$, $g\in S$.}
\label{rrgen}
\begin{tabular}{|l|l|l|l|l|l|}
\hline
\multirow{2}{*}{I. Structure} & \multirow{2}{*}{II. \texttt{GAP} ID} &III. Order$\vphantom{S^{S^S}}$  &   \multicolumn{2}{c|}{IV. Orbit sizes}     \\
\cline{4-5}& & \multicolumn{1}{|c|}{factored} &  $S\cong O_8^+(2)\vphantom{S^{S^S}g_{g_{g_g}}}$  &  $S\cong O_8^+(3)$ \\
\hline
  $3\vphantom{S^{S^S}}$                                     & [ 3, 1 ]      &  $3$         & $1$           &  $1$                 \\
  $3^2$                               & [ 9, 2 ]      &  $3^2$       & $56$          &  $728$               \\
  $A_4$                                     & [ 12, 3 ]     &  $2^2.3$     & $63$ ($\times 3$)     &  $351$ ($\times 3$)          \\
  $\mathrm{SL}_2(3)$                                 & [ 24, 3 ]     &  $2^3.3$     & $378$, $1512$ &  $44226$             \\
  $3^{1+2}_+$                               & [ 27, 3 ]     &  $3^3$       & $56$ ($\times 3$)     &  $17472$, $728$ ($\times 3$) \\
  $3 \times \mathrm{SL}_2(3)$                      & [ 72, 25 ]    &  $2^3.3^2$   & $1512$        &  $176904$            \\
  $(3^3\!:\! 2^2)\!:\!3$ & [ 324, 160 ]  &  $2^2.3^4$   & $2016$ ($\times 3$)   &  $78624$ ($\times 3$)        \\
  $[2^7]\!:\! 3^{1+2}_+$                       & [ 3456, n/a\tnote{\dag} \mbox{}\ ] &  $2^7.3^3$   & $1512$ ($\times 3$)   &  $176904$ ($\times 3$)       \\
  $[3^5]\vphantom{g_{g_{g_g}}}$ & [ 243, 3 ]    &  $3^5$       &    \multicolumn{1}{|c|}{---\ \tnote{\ddag}}           &  $157248$            \\
\hline 
\end{tabular}
\begin{tablenotes}
{\small 
\item[\dag] No \texttt{GAP} ID is available for groups of size $3456$.
\item[\ddag] Subgroups of order $243$ occur for $S\cong O_8^+(3)$ only.
}
\end{tablenotes}
\end{threeparttable}
\end{table}

As we can see from Column III, the order of each subgroup is a $\{2,3\}$-number. In particular, each subgroup is solvable and $\alpha(\rho)=3$ which proves $(i)$ for $S\cong O_8^+(2)$.

Then we show in \cite{RZGenTr} that $\varrho= (\rho x^3)^4 \in S\rho$ has order $3$ and is not conjugate to $\rho$ by calculating the order 
$$|C_S(\varrho)|=216=|\mathrm{PGU}_3(2)|$$ 
and using (\ref{crv}). Upon verifying that $\langle \varrho,\varrho^x\rangle=\langle S,\rho\rangle$, we obtain $\alpha(\varrho)\leqslant 2$. Clearly, we also have $\alpha(\varrho)\geqslant 2$, and so $\alpha(\varrho)=2$ thus proving $(ii)$ for $S\cong O_8^+(2)$.

Next, let $S\cong O_8^+(3)$. The extension $G={S\!:\!S_4}$ can be constructed in \texttt{GAP} \cite{GAP} using the function \texttt{AtlasGroup()} as a permutation group of degree $3360$. It has two generators $x_0$ and $y_0$ with orders $24$ and $20$, respectively, such that  $\langle x,y\rangle=S$, where  $x=x_0^4$, $y=y_0^4$, and $S=G'''$. We check that $\rho =(x_0y_0)^4$ has order $3$ and induces on $S$ a triality automorphism by calculating the order $$|C_S(\rho)|=4245696=|G_2(3)|$$ 
and using \cite{atlas}, cf. (\ref{crv}). 

The group $\langle \rho, \rho^x, \rho^y \rangle$ is the entire extension $\langle S, \rho\rangle$, which implies $\alpha(\rho)\leqslant 3$. As above, we calculate the orders of $2$-generated subgroups $\langle \rho,\rho^g\rangle$ for orbit representatives $\rho^g$ of $C_S(\rho)$ acting on the conjugacy class $\rho^S$. In this case, there are $18$ orbits resulting in subgroups of the same isomorphism types as for $O_8^+(2)$ (which is not surprising as ${O_8^+(2)\!:\!3}$ is a subgroup of ${O_8^+(3)\!:\!3}$, see~\cite{atlas}) plus one new type of order ${243=3^5}$. Again, since all the orders are $\{2,3\}$-numbers, we confirm $(i)$ for $S\cong O_8^+(3)$ as above.

Finally, we verify that $\varrho=(\rho x^4y)^4$ is a nontriality element of order~$3$ in the coset $S\rho$ by calculating the order 
$$|C_S(\varrho)|=5832=|[3^5].\mathrm{SL}_2(3)|$$ 
and using (\ref{crv}). The group $\langle \varrho, \varrho^x\rangle$ turns out to be the whole extension $\langle S,\rho\rangle$ and, therefore, $\alpha(\varrho)=2$ similarly to the case of $O_8^+(2)$. The claim in $(ii)$ holds.
\end{proof}

{\em Added on August 6, 2025.} After publishing the original version of this preprint, we learnt that C. Parker and J. Saunders had pointed out in \cite[Remark 5.2]{ParSau}
one of the two cases missing from the list of exceptions in \cite[Theorem A$\mbox{}^*$]{Guest}, namely, the triality automorphism of $O_8^+(3)$. They also indicated the existence of a gap in \cite[proof of Lemma 7]{Guest}. We are grateful to Prof. Parker for informing us about this.

\end{document}